\documentclass[a4paper]{amsart}

\usepackage{amscd,amsthm,amsfonts,amssymb,esint, amsmath}
\usepackage{tikz-cd}
\usepackage[colorlinks=true]{hyperref}
\usepackage{fullpage}
\usepackage[all]{xy}
\usepackage{color}
\usepackage{comment}
\usepackage{enumitem}
\usepackage[mathscr]{euscript}
\usepackage{mathtools}
\usepackage{mathrsfs}

    \newcommand{\BC}{{\mathbb {C}}}

     \newcommand{\BZ}{{\mathbb {Z}}}

    \newcommand{\CA}{{\mathcal {A}}} 
    \newcommand{\CC}{{\mathcal {C}}}

    \newcommand{\CO}{{\mathcal {O}}} 
     
    \newcommand{\CS}{{\mathcal {S}}} \newcommand{\CT}{{\mathcal {T}}}

     \newcommand{\RN}{{\mathrm {N}}}

    \newcommand{\RU}{{\mathrm {U}}}

    \newcommand{\fg}{{\mathfrak{g}}} 
     
     \newcommand{\fl}{{\mathfrak{l}}}

    \newcommand{\fu}{{\mathfrak{u}}}

    \newcommand{\End}{{\mathrm{End}}}

     \newcommand{\GL}{{\mathrm{GL}}}
    
    \newcommand{\Hom}{{\mathrm{Hom}}}

    \newcommand{\length}{{\mathrm{length}}}

    \renewcommand{\Re}{{\mathrm{Re}}}

    \newcommand{\st}{{\mathrm{st}}}

    \newcommand{\tr}{{\mathrm{tr}}}

    \newcommand{\wt}{\widetilde}
    \newcommand{\wh}{\widehat}

    \newcommand{\sk}{\medskip}
    \newcommand{\lra}{\longrightarrow}
    \newcommand{\ra}{\rightarrow} 
    \newcommand{\bs}{\backslash}
    
    \newcommand{\s}{\sk\noindent}

    \theoremstyle{plain}
    \newtheorem{thm}{Theorem}[section] 
    \newtheorem{lem}[thm]{Lemma}  \newtheorem{prop}[thm]{Proposition}

\theoremstyle{remark} 
\theoremstyle{exercise} 
\theoremstyle{example} 
\theoremstyle{question} 

\theoremstyle{assumption} 

    \numberwithin{equation}{section}

      \newcommand{\Quot}{{\mathrm{Quot}}}

    \renewcommand{\sp}{{\mathrm{sp}}}

\begin{document}

\title{Orbital Integrals on Unitary Groups and Zeta Values}

\author{Li Cai} 
\address[L.C.]{Academy for Multidisciplinary Studies, Beijing National Center for Applied Mathematics, Capital Normal University, Beijing, 100048, People's Republic of China
} 
\email{caili@cnu.edu.cn}

\author{Xuan Cui} 
\address[X.C.]{Academy for Multidisciplinary Studies, Beijing National Center for Applied Mathematics, Capital Normal University, Beijing, 100048, People's Republic of China
} 
\email{2250502173@cnu.edu.cn}

\begin{abstract}
We generalize Yun's interpolation formula for orbital integrals on 
$\GL_n$ in terms of Dedekind zeta functions to stable orbital integrals on unitary groups. Our approach combines Xiao's germ expansion of orbital integrals and the Jacquet-Rallis fundamental lemma. 
\end{abstract}

\subjclass[2020]{Primary: 11R42; Secondary: 22E35, 11S40}         

\keywords{
Stable orbital integrals; unitary groups; zeta functions of orders}

\maketitle

\tableofcontents

In \cite{Yun13}, Yun established a remarkable connection between orbital integrals on $\mathrm{GL}_n$ and  Dedekind zeta values. Roughly speaking, let $F$ be a number field and let $\gamma\in\mathfrak{gl}_n(F)$ be a regular semisimple element. Set $K=F[\gamma]$, which is a finite étale $F$-algebra. Then the orbital integrals associated with $\gamma$ admit an interpolation in terms of the Dedekind zeta function (Theorem \ref{thm-gl-global})
\[\zeta_K(s).\]

As observed by Arthur \cite{Art18}, this connection 
gives potential applications to the beyond endoscopy program of Langlands. More recently, this interpolation property of orbital integrals has played a crucial role in \cite{CD25}, where the leading term in the asymptotic counting of integral matrices with a prescribed characteristic polynomial is expressed in terms of  Dedekind zeta values.

It is natural to ask whether such an interpolation property extends to other reductive groups.  In this paper, we focus on 
(stable) orbital integrals on unitary groups.  Let $E/F$ be
a quadratic extension of number fields and $\eta$ the associated
quadratic character. Let $G$ be a
unitary group with respect to $E/F$ and $\fg$ its Lie algebra. 
Roughly speaking, for a regular semisimple element $\gamma \in \fg(F)$, we show that the stable orbital integral of $\gamma$ admits an interpolation in terms of the following quadratic $L$-function (Theorem \ref{thm-unitary-global})
\[L(s,\eta_K)\]
where $K = F[\gamma]$ and  $\eta_K = \eta \circ \RN_{K/F}$.

The key ingredient of our proof is a theorem of Xiao \cite{Xia19}
(Theorem \ref{thm-xiao-general}), giving the smooth matching for singular orbital integrals  in the relative trace formula approach to the Gan-Gross-Prasad conjecture for unitary groups.

\section{Orbital integrals and zeta values}\label{sec-orb}

\subsection{The linear group case}

We first recall some results of Yun \cite{Yun13}.

Let $F$ be a $p$-adic field. Denote by $G' = \GL_n$ with $\fg'$ its Lie algebra. We consider the conjugation action of $G'$ on $\fg'$,
\[x \cdot g' = (g')^{-1} x g', \quad x \in \fg', \quad g' \in G'.\] 
Denote by $\CS(\fg') = C_c^\infty(\fg')$ the space of Schwartz functions on $\fg'$.

For each regular semisimple $\gamma' \in \fg'$,
consider the orbital integral
\[\CO_{\gamma'}(f') = \int_{T_{\gamma'} \bs G'} f'(\gamma' \cdot g)dg, \quad f' \in \CS(\fg').\]
Here, $T_{\gamma'}$ is the stabilizer of $\gamma'$ in $G'$. The
measure is normalized such that the volumes of $G'(\CO)$  and the
maximal open compact subgroup of $T_{\gamma'}$ are both $1$. 

We now describe an interpolation property of this orbital integral in terms of certain zeta functions. We shall focus on the case $f' = 1_{\fg'(\CO)}$ with $1_{\fg'(\CO)}$ the characteristic function of $\fg'(\CO)$. Write $\CO_{\gamma'} = \CO_{\gamma'}(1_{\fg'(\CO)})$. 

Denote by $\chi(x) \in F[x]$ the characteristic polynomial of $\gamma'$. Then $\CO_{\gamma'} = 0$ unless $\chi(x) \in \CO[x]$. We assume this in the following. Denote by $K = F[x]/(\chi(x))$ and $R = \CO[x]/(\chi(x))$. In particular, $R$ is an order of $K$.

If $L_1,L_2$ are $\mathcal O$-lattices, denote by
\[\ell(L_1, L_2) = \mathrm{length}_{\mathcal O}(L_1/(L_1 \cap L_2))-\mathrm{length}_{\mathcal O}(L_2/(L_1 \cap L_2)), \quad 
    [L_1:L_2]=q^{\ell(L_1,L_2)}.\]
Denote by $R^\vee=\Hom_{\mathcal O_F}(R,\mathcal O_F)$. For each $j \geq 0$, let $\Quot_{R^\vee}^j$ be the set of fractional $R$-ideals $M \subset R^\vee$ with $\ell(R^\vee,M) = j$. 

Consider the following zeta function
\[J_R(s) = \sum_{j \geq 0} \# \Quot_{R^\vee}^j \ q^{-js}\]
and
\[\wt{J_R}(s) = \frac{J_R(s)}{[\CO_K:R]^{-s} \zeta_K(s)}.\]

The following is \cite[Theorem 2.5]{Yun13} and \cite[Equation (4.1)]{Yun13}. 
\begin{thm}[Local theory]\label{thm-gl-local} We have the following.
\begin{enumerate}
\item (Meromorphic continuation) 
			$\wt{J_R}(s) = [\CO_K:R]^sP_R(q^{-s})$, where $P_R(t) \in 1 + t\BZ[t]$ is of degree $2\length_\CO(\CO_K/R)$.
			In particular, if $R = \CO_K$, then $\wt{J_\gamma}(s) = 1$.
\item (Functional equation) $\wt{J_R}(s) = \wt{J_R}(1-s)$.
\item (Interpolation property) $\wt{J_R}(1) = \CO_\gamma$. 
\end{enumerate}
\end{thm}

In fact, the zeta function $J_R(s)$ can be expressed as a relative orbital integral considered by Jacquet. This viewpoint is useful for generalizing the above theory for $\GL_n$ to unitary groups.

 Denote by $V^*=F_n$  the space of row vectors. Consider the action of $G'$ on $\fg' \times V^*$ by 
\[(x,v) \cdot g = \left( g^{-1}x g, vg \right).\] 

For a regular semisimple $\gamma' \in \fg'$, choose a vector $v^* \in V^*$ cyclic for $\gamma'$ in the sense that  $v^*,v^*\gamma,\cdots,v^*\gamma^{n-1}$ gives a $F$-basis of $V^*$. Consider the following orbital integral
\[J_{\gamma'}(s,f') = \int_{G'} f'\left( (\gamma',v^*) \cdot g \right) |\det g|^s dg, \quad f' \in \CS(\fg' \times V^*).\]
 
In the following, we shall consider $f' = 1_{\fg'(\CO) \times L_0^*}$ with $L_0^* =\CO^n$ the standard lattice in $V^*$.
Denote by 
\[J_{\gamma'}(s) = [v^*R:L_0^*]^s J_{\gamma'}\left( s,1_{\fg(\CO) \times L_0^*} \right).\]

\begin{prop}[{\cite[Lemma 4.5]{Yun13}}] \label{prop-gl-jacquet}
    For $\Re(s) \gg 0$, the integral defining $J_\gamma(s)$
    is absolutely convergent. Moreover, we have
    \[J_{\gamma'}(s) = J_R(s).\]
    In particular, the definition of $J_{\gamma'}(s)$ is independent of the choice of the vector $v^* \in V^*$.
\end{prop}


Now, let $F$ be a number field. Let $\chi(x) \in \CO[x]$ be a separable polynomial. Denote by $K = F[x]/(\chi(x))$ and
$R = \CO[x]/(\chi(x))$, then $R$ is an order of $K$. Consider the zeta function for $R$
\[\zeta_R(s) = \sum_{M \subset R^\vee} \left( \# R^\vee/M \right)^{-s}\]
where $M$ runs over nonzero $R$-modules in  $R^\vee$ of finite index (which we call $R$-sublattices in $R^\vee$). Denote by
$\zeta_K(s)$ the Dedekind zeta function for $K$ (without the archimedean components).

The following is from the Eulerian product of $\zeta_R$  (\cite[(Equation (3.2)]{Yun13}).

\begin{thm}[Global theory]\label{thm-gl-global} Let $\gamma' \in \fg'$ with characteristic polynomial $\chi(x) \in \CO[x]$.  We have
	\[\zeta_R(s) = \prod_{v < \infty} J_R(s) = [\CO_K:R]^{-s}\zeta_K(s) \prod_{v < \infty} \wt{J_{R_v}}(s), \quad
	\frac{\zeta_R(s)}{\zeta_K(s)}\Big|_{s=1} = \frac{\CO_{\gamma'}(1_{\fg'(\wh{\CO})})}{[\CO_K:R]}.\]
    Here, the test function $1_{\fg'(\wh{\CO})} = \bigotimes_{v < \infty} 1_{\fg'(\CO_v)}$.
\end{thm}

\subsection{The unitary group case}\label{sec-unitary}

Let $F$ be a $p$-adic field. Let $E$ be a separable quadratic extension of $F$. Denote by $\,\bar{}\,$ the Galois involution with respect to
$E/F$. Let $\eta$ be the quadratic character associated to $E/F$.

Denote by $G$ the quasi-split unitary group associated with an $n$-dimensional
Hermitian space with respect to $E/F$
\[G = \left\{ g \in \GL_{n,E} \Big| gJ\bar{g}^t = J \right\}\]
where $J$ is the  anti-diagonal matrix with all entries $1$. Denote by $\fg$ the (twisted) Lie algebra of $G$. Then 
\[\fg = \left\{ x \in \fg\fl_{n,E} \Big| xJ= J\bar{x}^t  \right\}.\]
We consider the conjugate action of $G$ on $\fg$,
\[x \cdot g = g^{-1} x g, \quad x \in \fg, \quad g \in G.\] 
Denote by $\CS(\fg) = C_c^\infty(\fg)$ the space of Schwartz
functions on $\fg$.

For each regular semisimple $\gamma \in \fg$, consider the stable conjugacy class defined by 
$$\CC_{\gamma}=\{\delta\in \fg|\;\delta\; \text{is conjugacy to}\; \gamma \;\text{over}\; \bar F\}/\sim.$$
Here, $\delta_1\sim\delta_2$ if there exists $g\in G$ such that $\delta_1= \delta_2 \cdot g$.

Consider the stable orbital integral
\[\CO_\gamma^\st(f) = \sum_{\delta \in\mathcal C_ \gamma} \CO_{\delta}(f), \quad f \in \CS(\fg)\]
where for each $\delta \in \CC_\gamma$,
\[\CO_{\delta}(f) = \int_{T_\delta \bs G} f(\delta \cdot g)dg.\]
Here, $T_\delta$ is the stabilizer of $\delta$ in $G$. The
measure is normalized such that the volumes of $G(\CO)$  and the
maximal open compact subgroup of $T_{\gamma}$ are both $1$. 

In the following, we want to give an interpolating property of the stable orbital integral in terms of certain zeta functions.

We shall focus on the case where $E/F$ is unramified and $f = 1_{\fg(\CO)}$. Write $\CO_\gamma^\st = \CO_\gamma^\st(1_{\fg(\CO)})$. 
Denote by $\chi(x) \in F[x]$ the characteristic polynomial of $\gamma$.  Then $\CO_\gamma^\st = 0$ unless $\chi(x) \in \CO[x]$. We assume 
this in the following. Denote by $K = F[x]/(\chi(x))$ and $R = \CO[x]/(\chi(x))$, then $R$ is an order of $K$. Denote by 
\[\eta_K = \eta \circ \RN_{K/F}.\]
Then $\eta_K$ is the quadratic character associated to the quadratic extension $E \otimes_F K / K$.

Consider the following zeta function
\[J_R(s,\eta) = \sum_{j \geq 0} \# \Quot_{R^\vee}^j \eta(\varpi)^j q^{-js}\]
and
\[\wt{J_R}(s,\eta) = \frac{J_R(s,\eta)}{[\CO_K:R]^{-s} L(s,\eta_K)}.\]
Note that, in the above, we consider the $\eta$-twist for the series defining $J_R(s)$.

\begin{thm}[Local theory]\label{thm-unitary-local} We have the following.
\begin{enumerate}	
    \item (Meromorphic continuation) 
			$\wt{J_R}(s,\eta) = [\CO_K:R]^sP_R(\eta(\varpi)q^{-s})$, where $P_R(t) \in 1 + t\BZ[t]$ is the polynomial in Theorem \ref{thm-gl-local} (1). 
			In particular, if $R = \CO_K$, then $\wt{J_R}(s,\eta) = 1$.
    \item (Functional equation) $\wt{J_R}(s,\eta) = \wt{J_R}(1-s,\eta)$.
	\item (Interpolation property) $\wt{J_R}(1,\eta) = \CO_\gamma^\st$. 
\end{enumerate}
\end{thm}
\begin{proof}[The proof of (1) and (2)]
    
  If $E/F$ is split, we have $\eta = 1$ and this is exactly Theorem \ref{thm-gl-local} (1) and (2).

  Assume $E/F$ is an extension of fields. We claim that
  \[L(s,\eta_K) = \zeta_K(s \pm s_0)\]
  where $s_0 \in \BC$ such that $\eta(\varpi) = -1 = q^{s_0}$. 
  To see this, denote by $K = \prod_i K_i$, $\varpi_i$ a
  uniformizer of $K_i$, $\eta_{K_i} = \eta \circ \RN_{K_i/F}$ and $f_i$ the residue degree of $K_i/F$. Then
  \[\eta_{K_i}(\varpi_i) = (-1)^{f_i}.\]
  Therefore,
\[L(s,\eta_K)=\prod_{i}L(s,\eta_{K_i})=\prod_{i}\bigl(1-(-1)^{f_i}q^{-f_is}\bigr)^{-1}=\prod_{i}\bigl(1-q^{-f_i(s\pm s_0)}\bigr)^{-1}=\zeta_K(s\pm s_0).\]
   Now, we have
\[\wt{J_R}(s,\eta)=\frac{J_R(s,\eta)}{[\CO_K:R]^{-s} L(s,\eta_K)}
=\frac{\wt{J_R}(s\pm s_0)}{[\CO_K:R]^{\pm s_0}}.  \]    
    Applying Theorem \ref{thm-gl-local}, we obtain Theorem 
\ref{thm-unitary-local} (1) and (2).

\end{proof}

In the following, we consider the interpolation property 
of stable orbital integrals (Theorem \ref{thm-unitary-local} (3)).

Similar to Proposition \ref{prop-gl-jacquet}, we write $J_R(s,\eta)$ as certain relative orbital integral.

For a regular semisimple $\gamma' \in \fg'$, choose a vector $v^* \in V^*$ cyclic for $\gamma'$.  Consider the following orbital integral
\[J_{\gamma'}(s,f',\eta) = \int_{G'} f'\left( (\gamma',v^*) \cdot g \right) |\det g|^s \eta(\det g) dg, \quad f' \in \CS(\fg' \times V^*).\]

In the following, we shall consider $f' = 1_{\fg'(\CO) \times L_0^*}$.
Then $J_{\gamma'}(s,f',\eta)$ is nonvanishing only when
the characteristic polynomial $\chi(x)$ of $\gamma'$ belongs to $\CO[x]$. As before, denote by $K = F[x]/(\chi(x))$ and
$R = \CO[x]/(\chi(x))$. Denote by 
\[J_{\gamma'}(s,\eta) = \Delta(\gamma',v^*)[v^*R:L_0^*]^s J_{\gamma'}\left( s,1_{\fg'(\CO) \times L_0^*},\eta \right)\]
where
\[\Delta(\gamma',v^*) = \eta(\det(C(\gamma',v^*))), \quad C(\gamma',v^*) = \begin{pmatrix}
v^* \\
v^*\gamma' \\
\vdots \\
v^*\gamma'^{n-1}
\end{pmatrix}.\]

We also write
\[\wt{J_{\gamma'}}(s,\eta) = \frac{J_{\gamma'}(s,\eta)}{[\CO_K:R]^{-s} L(s,\eta_K)}.\]

\begin{prop}\label{prop-unitary-jacquet}  For $\Re(s) \gg 0$, the integral defining $J_{\gamma'}(s,\eta)$
    is absolutely convergent. Moreover, we have
    \[J_{\gamma'}(s,\eta) = J_R(s,\eta).\]
     In particular, the definition of $J_{\gamma'}(s,\eta)$ is independent on the choice of the vector $v^* \in V^*$.
\end{prop}
\begin{proof}
    The proof is analogous to that of Proposition \ref{prop-gl-jacquet}. For the convenience of the reader, we give a
    proof in Section \ref{sec-proof-unitary-jacquet}.
\end{proof}

The main ingredient in the proof of Theorem \ref{thm-unitary-local}(3) is the following result, which essentially follows from Xiao's Theorem 3.12 \cite{Xia19} and the Jacquet--Rallis fundamental lemma \cite{Yun11,BP21}.

\begin{thm}\label{thm-xiao} Let $\gamma \in \fg$ and $\gamma' \in \fg'$ be regular semisimple. Assume they are matching in the
	sense that they have the same characteristic polynomial. Then
	\[\CO_\gamma^\st = \wt{J_{\gamma'}}(0,\eta).\]
\end{thm}
\begin{proof}
    We shall give a proof in Section \ref{sec-proof-thm-xiao}.
\end{proof}

Now, Theorem \ref{thm-unitary-local} (3) follows from  Theorem \ref{thm-xiao}, Proposition \ref{prop-unitary-jacquet} and
the functional equation of $\wt{J_R}(s,\eta)$ (Theorem \ref{thm-unitary-local} (2)). 

Let $E/F$ be a quadratic extension of number fields. We assume that $E/F$ is unramified at every finite place. 

Let $\gamma \in \fg$ be regular semisimple. 
Let $\chi(x) \in \CO[x]$ be the characteristic separable polynomial of $\gamma$. Denote by $K = F[x]/(\chi(x))$ and
$R = \CO[x]/(\chi(x))$, then $R$ is an order of $K$. Consider the zeta function for $R$
\[J_R(s,\eta) = \sum_{M \subset R^\vee} \eta(M)\left( \# R^\vee/M \right)^{-s}.\]
Here, $M$ runs over nonzero $R$-sublattices in  $R^\vee$ and 
\[\eta(M) = \prod_{v < \infty} \eta_v(\varpi_v)^{j_v(M_v)}, \quad j_v(M_v) = \ell(R_v^\vee,M_v).\]

\begin{thm}[Global theory]\label{thm-unitary-global} Let $\gamma \in \fg$ with the characteristic polynomial $\chi(x) \in \CO[x]$.  We have
	\[J_R(s,\eta) = \prod_{v < \infty} J_{R_v}(s,\eta_v) = [\CO_K:R]^{-s}L^{(\infty)}(s,\eta_K) \prod_{v < \infty} \wt{J_{R_v}}(s,\eta_v), \quad
	\frac{J_R(s,\eta)}{L^{(\infty)}(s,\eta_K)}\Big|_{s=1} = \frac{\CO_\gamma^\st(1_{\fg(\wh{\CO})})}{[\CO_K:R]}.\]
    Here, $L^{(\infty)}(s,\eta_K) = \prod_{v < \infty} L(s,\eta_{K,v})$ is the finite part of the $L$-function
    $L(s,\eta_K)$ and the test function $1_{\fg(\wh{\CO})} = \bigotimes_{v < \infty} 1_{\fg(\CO_v)}$ with
    $\CO_\gamma^\st(1_{\fg(\wh{\CO})}) = \prod_{v < \infty}\CO_\gamma^\st(1_{\fg(\CO_v)})$.
\end{thm}
\begin{proof}
The proof is similar to that of Theorem \ref{thm-gl-global}. 

By the local-global correspondence for lattices, the map
\[M\mapsto (M_v)_{v<\infty}\]
induces a bijection between the set of $R$-sublattices in $R^\vee$  and the set of families of $R_v$-sublattices in $R_v^\vee$ such that $M_v=R_v^\vee$ for almost all $v$. 

For each $M \subset R^\vee$, we have
\[\#(R^\vee/M)=\prod_{v<\infty}\#(R_v^\vee/M_v)=\prod_{v<\infty}q_v^{j_v(M)}, \quad 
  \eta(M)=\prod_{v<\infty}\eta_v(\varpi_v)^{j_v(M)}.\]
Therefore, by the above local-global correspondence,
\[
\begin{aligned}
  J_R(s,\eta)
  &=\sum_{M\subset R^\vee}
      \eta(M)\bigl(\#(R^\vee/M)\bigr)^{-s} \\
  &=\prod_{v<\infty}
      \left(
        \sum_{M_v\subset R_v^\vee}
        \eta_v(\varpi_v)^{j_v(M_v)}
        q_v^{-s j_v(M_v)}
      \right)                                      \\
  &=\prod_{v<\infty}J_{R_v}(s,\eta_v).
\end{aligned}
\]

The remaining part of the Theorem follows from the
interpolate property of the stable orbital integral (Theorem
\ref{thm-unitary-local} (3)).

\end{proof}

\section{The proof}

Here, we give a proof of Proposition \ref{prop-unitary-jacquet} and Theorem \ref{thm-xiao}. Let $F$ be a $p$-adic field. Let $E$ be a separable quadratic extension of $F$. Denote by $\eta$ the character associated to $E/F$. When $E$ is a field, we assume $E/F$ is unramified. We shall keep notations in Section \ref{sec-orb}.

\subsection{The proof of Proposition \ref{prop-unitary-jacquet}}
\label{sec-proof-unitary-jacquet}

We shall only consider the case $\eta$ is non-trivial. Let
$\gamma' \in \fg'$ be regular semisimple with $K = F[\gamma']$.
Let $v^* \in V^*$ be cyclic to $\gamma'$.

Denote by $V=F^n$  the space of column vectors. There is a pairing between $V^*$ and $V$
\[
    \left\langle \cdot,\cdot\right\rangle:V^*\times V\to F,\qquad\left\langle v_1,v_2\right\rangle=v_1v_2.
\]

Let $L_0=\CO^n$ be the standard lattice in $V$.
\begin{lem}\label{lem-5}
We have the following bijection
\begin{align*}
    G'(F)/G'(\CO)&\stackrel{\sim}{\lra}\{\text{lattices in V}\}\\gG'(\CO)&\longmapsto \;\; L_g = gL_0.
\end{align*}
For $g\in G'(F)/G'(\CO)$, we have 
\begin{equation}
    g^{-1}\gamma'g\in\fg'(\CO) \;\;\text{if and only if }\;\; \gamma' L_g\subseteq L_g.
\end{equation}
Moreover, if $g$ satisfies (2.1), then
\begin{equation}
    v^*g\in L_0^*\;\;\text{if and only if }\;\; \left\langle v^*R,L_g\right\rangle \subseteq\CO.
\end{equation}

\end{lem}
\begin{proof}
First, we prove (2.1). Let $\delta=g^{-1}\gamma'g\in\fg'(\CO)$, equivalently we prove    
$\delta\in\fg'(\CO)$ iff $\delta L_0\subseteq L_0$.
If $\delta\in \fg'(\CO)$, then $\delta L_0\subseteq L_0$ is obvious. Conversely, choose $e_k=(0,...,1,...,0)^T\in L_0$, $1\leq k\leq n$. Then $\delta e_k\in L_0$ means every column of $\delta$ is in $L_0$, therefore $\delta\in \fg'(\CO)$.

 Then we prove (2.2). Note that (2.1) implies $RL_g\subseteq L_g$, hence \[\left\langle v^*R,L_g\right\rangle=\left\langle v^*,RL_g\right\rangle\subseteq\left\langle v^*,L_g\right\rangle=\left\langle v^*g,L_0\right\rangle.\]
By definition, $v^*g\in L_0^*$ is equivalent to $\left\langle v^*g,L_0\right\rangle \subseteq \CO$. Therefore $ \left\langle v^*R,L_g\right\rangle\subseteq\CO$ if and only if $v^*g\in L_0^*$.
\end{proof}

Since $v^*$ is cyclic to $\gamma'$, we can identify $V^*$ with $K$ by
\[\iota^*:K\stackrel{\sim}{\longrightarrow} V^*,\qquad\iota^*(a)=v^*a.\]
Consider the following pairing on $K$ 
\[(\cdot,\cdot): K\times K\to F,\qquad(x,y)= \tr_{K/F}(xy).\]

For any $v\in V$, there exists a unique $b\in K$ such that 
\[\langle\iota^*(a),v\rangle=(a,b), \quad a \in K.\]
This gives an isomorphism  
\[\iota:V\stackrel{\sim}{\longrightarrow} K,\qquad\iota(v)=b.\]
In particular, the following diagram is commutative
\[
\begin{tikzcd}
V^*  \times  V \arrow[d,"\iota",shift left=5] \arrow[r, "{\langle \cdot,\cdot\rangle}"] & F \arrow[d] \\
K \arrow[u, "\iota^*",shift left=4]  \times  K \arrow[r, "{(\cdot,\cdot)}"]                                & F          
\end{tikzcd}
\]

\begin{lem}\label{lem-6}
There is a bijection
    \begin{align*}
        \{\text{lattices $L_g$ in $V$}|\gamma' L_g\subseteq L_g,\left\langle v^*R,L_g\right\rangle \subseteq\CO\}&\stackrel{\sim}{\lra}\{\text{fractional $R$-ideal $M$ in $K$}|M\subseteq R^\vee\}\\L_g&\longmapsto\;\iota(L_g)
    \end{align*}
\end{lem} 
\begin{proof}Under the identification of $V^*\cong K\cong V$, we have $\left\langle v^*R,L_g\right\rangle =(R,\iota(L_g))$. Therefore, $\left\langle v^*R,L_g\right\rangle \subseteq\CO$ is equivalent to $(R,\iota(L_g))\subseteq\CO$. With the condition $\gamma' L_g\subseteq L_g$ which implies $RL_g\subseteq L_g$, $\iota(L_g)$ becomes a fractional $R$-ideal contained in $R^\vee$. The inverse map is $M\mapsto\iota^{-1}(M)$.
\end{proof}

\begin{lem}\label{lem-7}
    We have $\Delta(\gamma',v^*) =(-1)^{-\ell(v^*R,L_0^*)}$.
\end{lem}
\begin{proof}
    Since $v^*R$ is generated by the rows of the cyclic matrix $C(\gamma',v^*)$, we have
    \[\ell(v^*R,L_0^*)=-v(\det C(\gamma',v^*)).\]
\end{proof}

\begin{proof}[Proof of Proposition 1.5]

By definition
\[J_{\gamma'}(s,\eta) =
\Delta(\gamma',v^*)[v^*R:L_0^*]^s
\int_{G'}
\mathbf 1_{\fg'(\CO )}(g^{-1}\gamma'g)
\mathbf 1_{L_0^*}(v^*g)
\eta(\det g)|\det g|^s\,dg.\]
Note that the above integral is right $G'(\CO)$-invariant.
The Haar measure is normalized such that $G'(\CO)$ has volume $1$.
By  Lemma \ref{lem-5}, we have
\[\int_{G'}
\mathbf 1_{\fg'(\CO )}(g^{-1}\gamma'g)
\mathbf 1_{L_0^*}(v^*g)
\eta(\det g)|\det g|^s\,dg
= \sum_{\substack{\gamma' L_g\subseteq L_g \\\langle v^*R,L_g\rangle \subseteq\mathcal O}}\eta(\det g)|\det g|^s.\]

For each $L_g$, note that
\[\ell(R^\vee,\iota(L_g))=\ell(R^\vee,L_0)+\ell(L_0,\iota(L_g))=-\ell(e^*R,L_0^*)+v(\det g).\]
By Lemma \ref{lem-6}, we have
\[\sum_{\substack{\gamma' L_g\subseteq L_g \\\langle v^*R,L_g\rangle \subseteq\mathcal O}}\eta(\det g)|\det g|^s
= [v^*R:L_0^*]^{-s}(-1)^{\ell(v^*R,L_0^*)}\sum_{M\subseteq R^\vee}(-1)^{\ell(R^\vee,M)}[R^\vee:M]^{-s}.\]

By Lemma \ref{lem-7}, we obtain
\[J_{\gamma'}(s,\eta) = J_R(s,\eta).\]

The absolute convergence of $J_{\gamma'}(s,\eta)$ is given by Proposition \ref{prop-gl-jacquet}.

\end{proof}

\subsection{A consequence of Xiao's theorem}


Consider the following action of $G'$ on $\fg'\times V\times V^*$:
\[(x,v,v^*)\cdot g=(g^{-1}xg,g^{-1}v,v^* g).\]
An element $(x,v,v^*)\in \fg'\times V\times V^*$ is called regular semisimple if $\Delta(x,v,v^*)=\det\big(\langle v^*,x^{i+j}v\rangle_{0\leq i,j\leq n-1}\big)\neq0$. 
Here, we write $\big(\langle v^*,x^{i+j}v\rangle_{0\leq i,j\leq n-1}\big)$ as the matrix with $(i,j)$ entries $\langle v^*,x^{i+j}v\rangle$. 

For regular semisimple $(x,v,v^*) \in \fg' \times V \times V^*$, denote by \[\omega(x,v,v^*)=\eta(\det \{v,xv,...,x^{n-1}v\}).\]
Then $\omega((x,v,v^*) \cdot g) = \eta(\det g)\omega( (x,v,v^*))$.

Fix a regular semisimple element
$\gamma'\in\mathfrak{g}'.$
Let \[K=F[\gamma']=\prod_{i=1}^{m}K_i.\]
Consider the action of $K$ on $V$ by left multiplication. Let $e_i\in K$ be the idempotent associated with the $i$-th component $K_i$. We have the decomposition $V=\bigoplus_{i=1}^mV_i$, where $V_i=e_iV$. Similarly, we have $V^*=\bigoplus_{i=1}^m V^*_i$, where $V_i^*=V^*e_i$.

For each $1 \leq i \leq m$, denote by
\[\Delta_i(\gamma',v,v^*)=\Delta(\gamma_i',v_i,v^*_i)\]
where $(\gamma'_i,v_i,v^*_i)$ is the $i$-th component of 
$(\gamma',v,v^*)$ under the above decomposition.

Denote by $S_1=\{i|K_i\nsupseteq E\}$ and $S_2=\{i|K_i\supseteq E\}$. For each $\Lambda\subseteq S_1$ and given $v\in V$, denote by $v_{S_1\setminus\Lambda}\in V$ which equals $v$ in the $S_1\setminus\Lambda$ component and 0 otherwise. Similarly, given $v^*\in V^*$, denote by $v^*_{\Lambda}$ which equals $v^*$ in the $\Lambda$ component and 0 otherwise. Then $(\gamma',v_{S_1\setminus\Lambda},v_\Lambda^*) \in \fg' \times V \times V^*$ is not regular semisimple.

Let $T_{\gamma'}=K^{\times}=\prod_iK_i^{\times}$ be the stabilizer of $\gamma'$. Set $T_1=\prod_{i\in S_1}K_i^{\times}$. Consider the following (regularized) orbital integral
\begin{align*}
    \CO(f',(\gamma',v_{S_1\setminus\Lambda},v^*_{\Lambda}))=\int_{T_{\gamma'}\backslash G'} \left. \left( \int_{T_1} f((\gamma',v_{S_1\setminus\Lambda},v^*_{\Lambda})\cdot tg)\prod_{i\in S_1\setminus\Lambda}|t_i|^{-s}\prod_{i\in \Lambda}|t_i|^s\eta_K(t)dt\right)\right|_{s=0}\eta(\det g)dg,\\f'\in\mathcal S(\fg'\times V\times V^*).
\end{align*}
 The convergence of this integral is given by \cite[Lemma 2.4]{Xia19}.


Let $(W,\langle \cdot,\cdot\rangle)$ be a Hermitian space of rank $n$ for the extension $E/F$. Let $\RU(W)$ be the corresponding unitary group over $F$. Denote by 
\[\fu(W)=\{y\in\End _E(W)|\langle ya,b\rangle=\langle a,yb\rangle,\;a,b\in W\}\] 
the twisted Lie algebra of $\RU(W)$. Let $W_\sp$ be the quasi-split Hermitian space.

Consider the action of $\RU(W)$ on $\fu(W)\times W$
$$(y,w)\cdot g=(g^{-1}yg,g^{-1}w).$$
An element $(y,w)\in\fu(W)\times W$ is called regular semisimple if $\Delta(y,w)=\det\big(\langle w,y^{i+j}w\rangle_{0\leq i,j\leq n-1}\big)\neq0$. 

Fix a regular semisimple $\delta\in\mathfrak u(W)$. Let $$K=F[\delta]=\prod_{i=1}^{m}K_i.$$
Then $K$ can act on $W$ by left multiplication. Let $e_i\in K$ be the idempotent associated with the \(i\)-th component $K_i$. We have the decomposition $W=\bigoplus_{i=1}^m W_i$ of Hermitian spaces, where $W_i=e_iW$. 

For each $1 \leq i \leq m$, denote by
\[\Delta_i(\delta,w)=\Delta(\delta_i,w_i)\]
where $(\delta_i,w_i)$ is the $i$-th component of 
$(\delta,w)$ under the abve decomposition.


We give an explicit construction of conjugacy classes in a fixed
stable conjugacy class. Fix a regular semisimple element
$\gamma'\in\fg'$. Let
\[
K=F[\gamma']=\prod_{i=1}^{m}K_i,\qquad
E_i=E\otimes_F K_i
\]
and
\[\CA=\prod_{i=1}^{m}\big(K_i^\times/N_{E_i/K_i}(E_i^{\times})\big)
\cong
\prod_{i\in S_1}\mathbb Z/2\mathbb Z.\]
Denote by $\mathcal P^{S_1}$ the power set of $S_1$.  We have
the following pairing
\begin{align*}
    \langle \cdot, \cdot \rangle:\mathcal A\times\mathcal P^{S_1}&\to\{\pm1\},
    \qquad\langle\lambda,\Lambda\rangle=(-1)^{\sum_{i\in \Lambda}\epsilon_i}, \quad (\epsilon_i)_i \in \CA \cong \prod_{i \in S_1} \BZ/2\BZ.
\end{align*}

Consider
\[
\CC=\left\{(\delta,W)|\delta\in\fu(W),\
\delta\text{ is regular semisimple and }\chi_{\delta}(t)=\chi_{\gamma'}(t)
\right\}\big/\cong.
\]
Here, $\chi_{\delta}(t)$ and $\chi_{\gamma'}(t)$ are the characteristic polynomials of $\delta$ and $\gamma'$ respectively. Two pairs
$(\delta_1,W_1)$ and $(\delta_2,W_2)$ are equivalence if there exists an isomorphism $g:W_1 \cong W_2$  with $g\delta_1 g^{-1}=\delta_2$.

For $\lambda=(\lambda_1,\lambda_2,\ldots,\lambda_m)\in\prod_{i=1}^{m}K_i^\times$, consider the $E$-vector space $W_\lambda=\bigoplus_{i=1}^{m}E_i$
with the Hermitian form
\[
\left\langle
(a_1,\ldots,a_m),(b_1,\ldots,b_m)
\right\rangle_\lambda
=
\sum_{i=1}^{m}
\tr_{E_i/E}
\left(\lambda_i a_i\overline{b_i}\right).
\]
We denote by $\gamma'_\lambda$ the endomorphism of
$W_\lambda$ given by the multiplication by $\gamma'\in K$. Then $\gamma'_{\lambda}\in\fu(W_{\lambda})$.

\begin{prop}[{\cite[Proposition 3.6]{Xia19}}] \label{porp-1} 
    With the above notation, the map $\lambda\mapsto (\gamma'_{\lambda},W_\lambda)$
defines a bijection between $\CA$ and $\CC$. 
\end{prop}

The  main result in Xiao's thesis \cite[Theorem 3.12]{Xia19}   relates  (relative) singular orbital integrals  between  general linear groups and unitary groups. By applying the Jacquet--Rallis fundamental lemma to this relation, we obtain the following result.

\begin{thm}\label{thm-xiao-general}
    Fix a regular semisimple element
$\gamma'\in\fg'$. Let $(\gamma'_\mathbf 1,W_{\mathbf 1})$ be the conjugacy class corresponding to the identity in $\CA$ under Proposition \ref{porp-1}. Choose $w\in W_{\mathbf 1}$ cyclic for $\gamma'_{\mathbf 1}$. Moreover, choose $v$ and $v^*$ such that $\Delta_i(\gamma',v,v^*)=\Delta_i(\gamma'_{\mathbf 1},w)$ for $1 \leq i \leq m$. Then, for each $\Lambda\subseteq S_1$, we have
\begin{align}\label{eq-xiao}
    \omega(\gamma',v,v^*)\,2^{|S_1|}
\CO \big(1_{\fg'(\CO)\times L_0\times L_0^*},(\gamma',v_{S_1\setminus\Lambda},v^*_\Lambda)\big)
=
\sum_{\substack{\lambda\in\mathcal A\\W_\lambda\cong W_\sp}}
\langle\lambda,\Lambda\rangle
\CO_{\gamma'_{\lambda}}(1_{\fg(\CO)}).
\end{align}

\end{thm}
\begin{proof}
    In \cite[Theorem 3.12]{Xia19}, (\ref{eq-xiao}) holds
    for general matching functions $f' \sim (f_W)_W$ with
    $f' \in \CS(\fg' \times V \times V^*)$ and each $f_W \in \CS(\fu(W) \times W)$ as following
    \[ \omega(\gamma',v,v^*)\,2^{|S_1|}
\CO \big(f',(\gamma',v_{S_1\setminus\Lambda},v^*_\Lambda)\big)
=
\sum_{\lambda\in \CA}
\langle\lambda,\Lambda\rangle
\CO_{\gamma'_{\lambda}}(f_{W_\lambda})\]
with
\[\CO_{\gamma'_{\lambda}}(f_{W_\lambda}) = 
\int_{T_{\gamma'_\lambda} \bs \RU(W_\lambda)} f_{W_\lambda}(\gamma'_\lambda \cdot g,0) dg.\]
    
By the Jacquet-Rallis fundamental lemma (\cite[Theorem 2.10]{Xia19}), 
$1_{\fg'(\CO)\times L_0\times L_0^*}$  matches with
$(1_{\fg(\CO) \times \CO_E^n},0)$. We obtain the result.
\end{proof}

We give several remarks:
\begin{enumerate}
    \item The existence of $v$ and $v^*$ will be given in Lemma \ref{lem-2}.
    \item The condition that $w \in W_{\bf 1}$ is cyclic for $\gamma_{\bf 1}'$
    implies that $\Delta_i(\gamma_{\bf 1}',w) \not=0$ for all $1 \leq i \leq m$. Moreover, this implies $(\gamma_{\bf 1}',w)$ and $(\gamma',v,v^*)$ are both regular semisimple.
    \item For each $1 \leq i \leq m$, $\eta(\Delta_i(\gamma_{\bf 1}',w))$ is independent on the choice of the cyclic vector $w \in W_{\bf 1}$.
    \item The left hand side of (\ref{eq-xiao}) depends on
    $v$ and $v^*$ only through $\eta(\Delta_i(\gamma',v,v^*))$, $1 \leq i \leq m$. Hence, by (3),  the left hand side
     of (\ref{eq-xiao}) is independent on the choice of $w$, $v$ and $v^*$.
\end{enumerate}

\subsection{The proof of Theorem \ref{thm-xiao}}
\label{sec-proof-thm-xiao}

Let $W$ be a Hermitian space of rank $n$ with respect to $E/F$.
Denote by $\det W \in F^\times/\RN_{E/F} E^\times$ the determinant
of $W$. For two such Hermitian spaces $W_1$ and $W_2$, $W_1$ is
isomorphic to $W_2$ if and only if $\eta(\det W_1) = \eta(\det W_2)$. We have $\eta(\det W_\sp)) = 1$.

\begin{lem}\label{lem-8}
   Let $W_\lambda$ be Hermitian space associated to $\lambda \in \CA$ given before Proposition \ref{porp-1}. Then
   \[\eta(\det W_\lambda) = \eta(\det W_{\bf 1}) \prod_{i \in S_1}
   \eta_{K_i} (\lambda_i).\]
\end{lem}
\begin{proof}
For a Hermitian space $(W,\langle \cdot,\cdot \rangle)$ of
rank $n$ with respect to $E/F$, 
\[\det W = \det H \in F^\times/\RN_{E/F} E^\times\]
where the Gram matrix $H = (\langle w_i,w_j \rangle)_{i,j}$ is
given with a basis $w_1,\ldots, w_n$ of $W$. 

Let $\lambda \in \CA$ be represented by $(\lambda_i)_i \in \prod_{i=1}^m K_i^\times$. Then 
\[(W_\lambda, \langle \cdot,\cdot \rangle_\lambda) = 
\left( \bigoplus_{i=1}^m E_i, (\langle \cdot,\cdot \rangle_{\lambda_i})_i \right)\]
with each
\[\langle a,b \rangle_{\lambda_i} = \tr_{E_i/E}(\lambda_i a\bar{b}), \quad a,b \in E_i.\]

We fix a $E$-basis for each $E_i$, $1 \leq i \leq m$. 
For any $\lambda \in \CA$, let $H_\lambda$ and $H_{\lambda,i}$, $1 \leq i \leq m$ be the Gram matrices of $E$ and $E_i$ with respect to this basis. Then
\[\det H_\lambda = \prod_{i=1}^m \det H_{\lambda,i}.\]

Moreover, we have
\[\det H_{\lambda,i} = \RN_{K_i/F}(\lambda_i) \det H_{{\bf 1},i}.\]

It follows that
\[
\det W_\lambda=\det W_{\mathbf 1}\prod_{i=1}^m\RN_{K_i/F}(\lambda_i),
\quad 
\eta(\det W_\lambda)=\eta(\det W_{\mathbf 1})
\prod_{i\in S_1}\eta_{K_i}(\lambda_i).
\]

\end{proof}

\begin{lem}\label{lem-1} 
We keep notations as in Theorem \ref{thm-xiao-general}. 
Let $\gamma \in \fg$ be regular semisimple which is conjugate
to $\gamma' \in \fg'$ over $\bar{F}$. Then 
\begin{align}\label{equ-lem-1}
    \sum_{\substack{\lambda\in\mathcal A\\W_\lambda\cong W_\sp}}\langle\lambda,S_1\rangle \CO_{\gamma'_{\lambda}}( 1_{\fg(\CO)})=\epsilon(W_{\mathbf 1})\CO_{\gamma}^{\st}.
\end{align}
Here, 
\[\epsilon(W)=\begin{cases}1,&\text{if} \;W\cong W_\sp \\-1,&\text{if} \;W\ncong W_\sp \end{cases}.\]
\end{lem}
\begin{proof}

By Lemma \ref{lem-8}, for each $\lambda \in \CA$, we have
\[\langle\lambda,S_1\rangle=(-1)^{\sum_{i\in S_1}\lambda_i}=\prod_{i\in S_1}\eta_{K_i}(\lambda_i) = \eta(\det W_\lambda)\eta(\det W_{\bf 1})^{-1}.\]

Take $\Lambda = S_1$ in Theorem \ref{thm-xiao-general}. 
The left side of (\ref{equ-lem-1}) becomes
\[\sum_{\substack{\lambda\in\CA\\W_\lambda\cong W_\sp}}
\langle\lambda,S_1\rangle \CO_{\gamma'_{\lambda}}( 1_{\fg(\CO)})=\sum_{\substack{\lambda\in\mathcal A\\W_\lambda\cong W_\sp}}
\eta(\det W_\sp)\eta(\det W_{\mathbf 1})^{-1}\CO_{\gamma'_{\lambda}}(1_{\fg(\CO)})=\epsilon(W_{\mathbf 1})\sum_{\substack{\lambda\in\mathcal A\\W_\lambda\cong W_\sp}}
\CO_{\gamma'_{\lambda}}( 1_{\fg(\CO)}).\]

By the bijection between $\CA$ and $\CC$ (Proposition \ref{porp-1}),
we have
\[
\sum_{\delta \in \CC_\gamma} \CO_\delta(1_{\fg(\CO)})
 = \CO_\gamma^\st.\]
 The Lemma holds.
\end{proof}
\begin{lem}\label{lem-2}
Let $\gamma'\in\fg'$ be regular semisimple. Choose $w\in W_{\mathbf 1}$ cyclic for $\gamma_{\mathbf 1}'$. Let $v^*\in V^*$ be cyclic for $\gamma'$. Then there exists a vector $v\in V$ such that $(\gamma',v,v^*)$ satisfies  $\Delta_i(\gamma',v,v^*)=\Delta
_i(\gamma_{\mathbf 1}',w)$, $1\leq i\leq m$ and $\omega(\gamma',v,v^*)=\epsilon(W_\mathbf{1})\Delta(\gamma',v^*)$.
\end{lem}
\begin{proof}

Since $v^* \in V^*$ is cyclic for $\gamma'$, 
$\{v^*,v^*\gamma',\ldots,v^*(\gamma')^{n-1}\}$ is a basis of $V^*$. There exists a unique $v \in V$
such that
\[\langle v^*,(\gamma')^jv \rangle = \langle w, (\gamma_1')^jw \rangle, \quad 0 \leq j \leq n-1.\]
Moreover, as $\gamma'$ and $\gamma'_1$ have the same characteristic
polynomial,
\begin{align}\label{eq-lem-2}
\langle v^*,P(\gamma')v \rangle = \langle w, P(\gamma_1')w \rangle, \quad P(t) \in F[t].
\end{align}

For each $1 \leq i \leq m$, denote by $e_i \in K = \prod_{i=1}^m K_i$  the idempotent associated to $K_i$. We have
\[\Delta_i(\gamma',v,v^*) = \det\left( \langle v_i^*, (\gamma_i')^{\ell+j} v_i \rangle_{\ell,j}\right)
= \det\left( \langle v^*, e_i(\gamma')^{\ell+j} v \rangle_{\ell,j}\right).
\]
Similarly, we have
\[\Delta_i(\gamma_{\bf 1}',w) = \det\left( \langle w, e_i(\gamma_{\bf 1})^{\ell+j}w \rangle_{\ell,j}\right).\]
By (\ref{eq-lem-2}), we obtain that
\[\Delta_i(\gamma',v,v^*) = \Delta_i(\gamma_{\bf 1}',w).\]

Denote by
\[H = (\langle v^*, (\gamma')^{i+j} v \rangle)_{i,j} =  C(\gamma',v^*)B, \quad B = (v,\gamma'v,\ldots,(\gamma')^{n-1}v).\]
Then
\[\eta(\det H) = \Delta(\gamma',v^*)\omega(\gamma',v,v^*).\]
On the other hand,
\[H= (\langle w, (\gamma_{\bf 1}')^{i+j} w \rangle)_{i,j} = 
(\langle (\gamma_{\bf 1}')^iw, (\gamma_{\bf 1}')^{j} w \rangle)_{i,j}\]
is the Gram matrix for $W_{\bf 1}$ with respect to the basis
$\{w,\gamma_{\bf 1}'w,\ldots, (\gamma_{\bf 1}')^{n-1}w \}$.
In particular,
\[\eta(\det H) = \eta(\det W_{\bf 1}) = \epsilon(W_{\bf 1}).\]
We obtain the Lemma.

\end{proof}

Let $\gamma' \in \fg'$ be regular semisimple with $K = F[\gamma'] = 
\prod_i K_i$. Denote by 
\[S_1 = \{i | K_i \not\supset E\}, \quad S_2 = \{i | K_i \supset E\}.\]
Denote by $T_{\gamma'}$ the stabilizer of $\gamma'$ in $G'$ and 
$T_1 = \prod_{i \in S_1} K_i^\times \subset T_{\gamma'}$. Let $v^* \in V^*$ be cyclic for $\gamma'$.

\begin{lem}\label{lem-3}
For any \(\phi\in\CS(V^*)\), we have
\begin{align}\label{equ-lem-3}
\left. 
\frac{\int_{T_{\gamma'}}\phi(v^*t)\eta_K(t)|t|_K^s\,dt}{\prod_{i\in S_2} \zeta_{K_i}(s)}
\right|_{s=0}
=
\left.
\int_{T_1}\phi(v^*_{S_1}t)\eta_K(t)
|t|_K^s\,dt
\right|_{s=0}.
\end{align}
\end{lem}
\begin{proof}
   We identify $K$ with $V^*$ by $t \mapsto v^* t$. We then may
   assume $\phi = \otimes_i \phi_i$ with $\phi_i \in \CS(K_i)$. 
   In particular,
   \[\int_{T_{\gamma'}}\phi(v^*t)\eta_K(t)|t|_K^s\,dt = 
   \prod_i \int_{K_i^\times} \phi_i(t)\eta_{K_i}(t)|t|_{K_i}^s dt.\]
   The character $\eta_{K_i} = 1$ if and only if $i \in S_2$. 
   For $i \in S_2$, we have
   \[\left. \frac{\int_{K_i^\times} \phi_i(t)|t|_{K_i}^s dt}
   {\zeta_{K_i}(s)} \right|_{s=0} = \phi_i(0), \quad \phi_i \in \CS(K_i).\]
   Therefore,
   \[\left. 
\frac{\int_{T_{\gamma'}}\phi(v^*t)\eta_K(t)|t|_K^s\,dt}{\prod_{i\in S_2} \zeta_{K_i}(s)}
\right|_{s=0}
= \prod_{i \in S_2} \phi_i(0) \left. \left(\prod_{i \in S_1} 
\int_{K_i^\times} \phi_i(t)\eta_{K_i}(t)|t|_{K_i}^s dt\right)\right|_{s=0} = \left.
\int_{T_1}\phi(v^*_{S_1}t)\eta_K(t)
|t|_K^s\,dt
\right|_{s=0}.\]

\end{proof}

\begin{lem}\label{lem-4}
    We have
    \[\wt{J_{\gamma'}}(0,\eta)=\Delta(\gamma',v^*)2^{|S_1|}\CO\big(\mathbf{1}_{\fg'(\CO)\times L_0\times L_0^*},(\gamma',0,v^*_{S_1})\big).\]
\end{lem}
\begin{proof}
   By definition,
   \[\wt{J_{\gamma'}}(s,\eta) = \frac{J_{\gamma'}(s,\eta)}{
   [O_K: R]^{-s}L(s,\eta_K)}\]
   with
   \[J_{\gamma'}(s,\eta) = \Delta(\gamma',v^*)[v^*R:L_0^*]^s
   J_{\gamma'}(s,1_{\fg'(\CO) \times L_0^*},\eta).\]

   We have
   \[J_{\gamma'}(s,1_{\fg'(\CO) \times L_0^*},\eta)
   = \int_{T_{\gamma'} \bs G'} 1_{\fg'(\CO)}(\gamma' g)|\det g|^s
   \eta(\det g) \left(\int_{T_{\gamma'}} 1_{L_0^*}(v^* tg)\eta_K(t)|t|_K^sdt\right) dg.\]
   By Lemma \ref{lem-3}, we have
   \[
   \begin{aligned}
   \left. \frac{J_{\gamma'}(s,1_{\fg'(\CO) \times L_0^*},\eta)}
   {L(s,\eta_K)} \right|_{s=0} &= 2^{|S_1|} 
   \int_{T_{\gamma'} \bs G'} 1_{\fg'(\CO)}(\gamma' g)|\det g|^s
   \eta(\det g) \left.\left(\frac{\int_{T_{\gamma'}} 1_{L_0^*}(v^* tg)\eta_K(t)|t|_K^sdt}{\prod_{i \in S_2} \zeta_{K_i}(s)}\right)
   \right|_{s=0}dg \\
   &=2^{|S_1|} 
   \int_{T_{\gamma'} \bs G'} 1_{\fg'(\CO)}(\gamma' g)|\det g|^s
   \eta(\det g) \left.\left( \int_{T_1} 1_{L_0^*}(v^*_{S_1}tg)
   \eta_K(t)|t|_K^s dt\right)
   \right|_{s=0}dg \\
   &= 2^{|S_1|} \CO(1_{\fg'(\CO) \times L_0 \times L_0^*},(\gamma',0,v_{S_1}^*)).
   \end{aligned}
   \]
   We obtain the Lemma.
\end{proof}

\begin{proof}[The Proof of Theorem \ref{thm-xiao}]
Using Lemma \ref{lem-2} and Lemma \ref{lem-4}, we have
\begin{align*}
   \omega(\gamma',v,v^*)2^{|S_1|} \CO\big(\mathbf{1}_{\fg'(\CO)\times L_0\times L_0^*},(\gamma',0,v^*_{S_1})\big)=\epsilon(W_{\mathbf 1})\wt{J_{\gamma'}}(0,\eta).
\end{align*}
On the other hand, by Lemma \ref{lem-1}, we have
\[\sum_{\substack{\lambda\in\CA\\W_\lambda\cong W_\sp}}\langle\lambda,S_1\rangle \CO_{\gamma'_\lambda}(\mathbf{1}_{\fg(\CO)})=\epsilon(W_{\mathbf 1})\CO_{\gamma}^{\st}.\]
Applying Theorem \ref{thm-xiao-general} with $\Lambda = S_1$, we obtain
\[\wt{J_{\gamma'}}(0,\eta)=\CO_\gamma^{\st}.
\]
\end{proof}

\section{A counting problem}

Here, we give a brief discussion of a potential application of our results to the counting problem of integral points.

Let $E/F$ be a quadratic extension of number fields with the quadratic character $\eta$. Denote by $\,\bar{}\,$ the Galois involution for $E/F$. 

Let $J$ be the anti-diagonal matrix with all entries 1. 
As in Section \ref{sec-unitary}, we consider the unitary group
\[G = \left\{ g \in \GL_{n,E} \Big| gJ\bar{g}^t = J \right\}.\]
and the (twisted) Lie algebra 
\[\fg = \left\{ x \in \fg\fl_{n,E} \Big| xJ= J\bar{x}^t  \right\}.\]

Let $\gamma \in \fg$ be regular semisimple. Denote by $\chi(t)$ the characteristic polynomial of $\gamma$. Then $\chi(t) \in F[t]$. 
Consider the affine variety
\[X = \left\{ x \in \fg \Big| \det(tI_n - x) = \chi(t) \right\}.\]

The unitary group $G$ acts on $X$ by conjugation
\[x \cdot g  = g^{-1}xg.\]
Denote by $T = T_\gamma$ the stabilizer of $\gamma$ in $G$. Then $T$ is a torus in $G$. As varieties over $F$, $X = T \bs G$. 

Assume $\gamma \in \fg(\CO)$ so that $\chi(t) \in \CO[t]$. Denote
by $K = F[t]/(\chi(t))$ and $R = \CO[t]/(\chi(t))$. Write
$\eta_K = \eta \circ \RN_{K/F}$.

For each positive number $\CT$, consider the counting function
\[N(X,\CT) = \left\{ x \in X(\CO) \Big| ||x|| \leq \CT \right\}\]
where $||\cdot||$ is fixed archimedean norm on $X$. 

The ergodic theory can be applied to study the asymptotic expansion
of $N(X,\CT)$ as $\CT \ra \infty$. Denote by $C(\chi)$ the leading coefficient for this asymptotic expansion. We want to express the constant $C(\chi)$ in terms of certain Dedekind zeta value.  

Following the method in  \cite{CD25} for the general linear group case, $C(\chi)$ is related to stable orbital integrals on unitary groups. 

By the above generalization of Yun's work to unitary groups (Theorem \ref{thm-unitary-global}), for the case $E/F$ is unramified everywhere, $C(\chi)$ is related to certain special value of
\[ L(s,\eta_K)\]
or more precisely
\[\left. \frac{J_R(s,\eta)}{L^{(\infty)}(s,\eta_K)} \right|_{s=1}.\]

\medskip
\noindent\textbf{Acknowledgments.}
The authors thank Taiwang Deng for helpful discussions
related to this work. L. Cai is supported by the National Key R\&D Program of China No. 2023YFA1009702 and National Natural Science Foundation of China, No. 12371012. 

\s{\bf Use of artificial intelligence.} ChatGPT was used to assist with language editing and checking mathematical arguments.

\end{document}